\documentclass[12pt,leqno]{amsart}  
\usepackage{amsmath,amstext,amsthm,amssymb,amsrefs}
\usepackage[colorlinks=true]{hyperref}
\usepackage{pgf,pgfarrows,pgfshade}
\usepackage{caption}   
\numberwithin{equation}{section}

\allowdisplaybreaks
\swapnumbers
\theoremstyle{plain}

\def\seq{\lesssim }

\def\bthm#1.#2 #3\ethm{
\begin{\ifx#1ttheorem\fi%
\ifx#1llemma\fi%
\ifx#1ccorollary\fi%
\ifx#1pproposition\fi%
\ifx#1ddefinition\fi}
\label{#1.#2}  
#3 \end{\ifx#1ttheorem\fi%
\ifx#1llemma\fi%
\ifx#1ccorollary\fi%
\ifx#1pproposition\fi%
\ifx#1ddefinition\fi}}

\def\t#1/{Theorem~\ref{t#1}}    
\def\c#1/{Corollary~\ref{c#1}}   
\def\l#1/{Lemma~\ref{l#1}}         
\def\s#1/{Section~\ref{s#1}}      
\def\e#1/{(\ref{e#1})}
\def\d#1/{Definition~\ref{d#1}}
\def\f#1/{Figure~\ref{f#1}}

\def\Label #1 {\label{#1}}

\def\norm#1.#2.{\lVert#1\rVert_{#2}}
\def\Norm#1.#2.{\bigl\lVert#1\bigr\rVert_{#2}}
\def\NOrm#1.#2.{\Bigl\lVert#1\Bigr\rVert_{#2}}
\def\NORm#1.#2.{\biggl\lVert#1\biggr\rVert_{#2}}
\def\NORM#1.#2.{\Biggl\lVert#1\Biggr\rVert_{#2}}

\def\ip#1,#2.{\langle #1,#2\rangle}
\def\Ip#1,#2.{\bigl\langle#1,#2\bigr\rangle}
\def\IP#1,#2.{\Bigl\langle#1,#2\Bigr\rangle}

\def\abs#1{\lvert#1\rvert}
\def\Abs#1{\bigl\lvert#1\bigr\rvert}
\def\ABs#1{\Bigl\lvert#1\Bigr\rvert}

\newcommand{\zc}{\ensuremath{\psi}}

\newcommand{\zd}{\ensuremath{\delta}}
\newcommand{\zD}{\ensuremath{\Delta}}
\newcommand{\ze}{\ensuremath{\epsilon}}
\newcommand{\zve}{\ensuremath{\varepsilon}}
\newcommand{\zvf}{\ensuremath{\varphi}}
\newcommand{\zf}{\ensuremath{\phi}}
\newcommand{\zF}{\ensuremath{\Phi}}

\newcommand{\zI}{\ensuremath{\infty}}
\newcommand{\zl}{\ensuremath{\lambda}}
\newcommand{\zk}{\ensuremath{\kappa}}

\newcommand{\zt}{\ensuremath{\tau}}
\newcommand{\zw}{\ensuremath{\omega}}
\newcommand{\zW}{\ensuremath{\Omega}}
\newcommand{\zx}{\ensuremath{\xi}}

\newcommand{\zs}{\ensuremath{\sigma}}

\newcommand{\zr}{\ensuremath{\rho}}

\newcommand{\zp}{\ensuremath{\pi}}
\newcommand{\zq}{\ensuremath{\chi}}
\newcommand{\zz}{\ensuremath{\zeta}}

\def\z#1#2{\ifcase#1 {{\mathcal {#2}}}  
\or {{\tilde{#2}}}                    
\or { {\boldsymbol{#2}}}                  
\or{{\widetilde{#2}}}                   
\or {{\acute{#2}}}
\or {{\grave{#2}}}
\or {{\bar{#2}}}
\or {\dot{#2}}
\or {\overline{#2}}
\or {\underline{#2}}\fi}
 
\def\ZR{\ensuremath{\mathbb R}}
\def\ZZ{\ensuremath{\mathbb Z}}

\def\ZN{\ensuremath{\mathbb N}}
\def\ZT{\ensuremath{\mathbb T}}

\def\vZT{\ensuremath{\mathcal T}}

\def\vZT{\ensuremath{\mathcal T}}

\def\Xint#1{\mathchoice
   {\XXint\displaystyle\textstyle{#1}}%
   {\XXint\textstyle\scriptstyle{#1}}%
   {\XXint\scriptstyle\scriptscriptstyle{#1}}%
   {\XXint\scriptscriptstyle\scriptscriptstyle{#1}}%
   \!\int}
\def\XXint#1#2#3{{\setbox0=\hbox{$#1{#2#3}{\int}$}
     \vcenter{\hbox{$#2#3$}}\kern-.5\wd0}}

\def\dashint{\Xint-}

 \def\inr{\int_{\ZR}}        
 
\def\exp#1 {{\text{\rm e}}^{#1}}  
\def\bw#1{{\boldsymbol \omega}_{s#1}}  

\def\tree#1{$#1$--tree} 
\def\T{ \ensuremath{{\mathbf T}}  }
\def\S { \ensuremath{{\mathbf S}} }
\def\ipf{ \ip f,\zvf_{s}.}

\def\mid{\,:\,}

\def\md#1#2\emd{\ifx0#1
\begin{equation*} #2 \end{equation*}\fi  
\ifx1#1\begin{equation}#2\end{equation}\fi   
\ifx2#1\begin{align*}#2\end{align*}\fi   
\ifx3#1\begin{align}#2\end{align}\fi    
\ifx4#1\begin{gather*}#2\end{gather*}\fi  
\ifx5#1\begin{gather}#2\end{gather}\fi   
\ifx6#1\begin{multline*}#2\end{multline*}\fi  
\ifx7#1\begin{multline}#2\end{mutline}\fi  
}

\def\prm#1{{\operatorname{\rm ann}}(#1)}
\def\ind#1{ {\mathbf 1}_{#1}}
\def\size#1{\operatorname{\rm size}(#1)}

\def\sh#1{\operatorname{\rm sh}(#1)}

\def\dense#1{\operatorname{\rm dense}(#1)}

\def\sh#1{\operatorname{sh}(#1)}  
\def\scl#1{\operatorname{\rm scl}(#1)}      

\begin{document}

\title[Maximal Directional Hilbert Transform]
{Maximal Theorems for the Directional Hilbert Transform on the Plane}
 \author{Michael T. Lacey}

\thanks{The research  of both authors was supported in part by  NSF grants; one of use (M.L.) is also supported by the 
 Guggenheim Foundation.  Some of this research was completed during 
research stays by M.L.~at the Universite d'Paris-Sud, Orsay, and  the  Erwin Schr\"odinger 
Institute of Vienna Austria.  The generosity of both is gratefully acknowledged.
}

\address{Michael Lacey\\
School of Mathematics\\
Georgia Institute of Technology\\
Atlanta,  GA 30332}

\email{lacey@math.gatech.edu}

\author{Xiaochun Li}

\address{Xiaochun Li\\
Department of Mathematics\\
University of California, Los Angeles\\
Los Angeles CA }

\email{xcli@math.ucla.edu}

\curraddr{Xiaochun Li\\
School of Mathematics\\
Institute for Advanced Study\\
Princeton, NJ  08540}

\date{\today}


\begin{abstract}{For a Schwartz function $f$ on the plane  and a non-zero $v\in\ZR^2$ define the Hilbert transform of $f$ in the 
direction $v$ to be 
\begin{equation*}
\operatorname  H_vf(x)=\text{p.v.}\int_\ZR f(x-vy)\; \frac{dy}y
\end{equation*}
Let $\zeta$ be a Schwartz function with frequency support in the annulus $1\le\abs \xi\le2$, and ${\boldsymbol \zeta}f=\zeta*f$. 
We prove that   the maximal operator  $\sup_{\abs v=1}\abs{ \operatorname H_v{\boldsymbol \zeta} f}$ maps $L^2$ into weak $L^2$, and $L^p$ into $L^p$ for $p>2$.  The $L^2$ estimate is sharp.  The method of proof is based upon techniques related to the pointwise convergence of Fourier series. 
Indeed, our main theorem implies this result on Fourier series.
} \end{abstract}

\maketitle

\section{Introduction, Principal Theorem}

Our interest is in the directional Hilbert transform applied in a choice of directions of the plane. Thus, 
for $v\in\ZR^2-\{0\}$, set 
\begin{equation*}
\operatorname H_vf(x)=\text{p.v.}\int_\ZR f(x-vy)\; \frac{dy}y
\end{equation*}
This definition is independent of the length of $v$, and below we shall only concern ourselves with $\abs v=1$. 
Let $\zeta$ be a Schwartz function with frequency support in the annulus $1\le\abs \xi\le2$, and ${\boldsymbol \zeta}f=\zeta*f$. 
Our Theorem is

\bthm t.main  The maximal operator  $\operatorname H^*_vf:=\sup_{\abs v=1}\abs{ \operatorname H_v {\boldsymbol \zeta}f}$ maps $L^2$ into weak $L^2$, and $L^p$ into $L^p$ for $p>2$. 
\ethm 

This theorem is a complement to a corresponding result, due to J.~Bourgain \cite{MR92g:42010},
for the directional maximal function, namely 
\begin{equation}  \label{e.dir-max}
\operatorname M^*f(x):=\sup_{\abs v=1} \sup_{t>0} (2t)^{-1}\int_{-t}^t \abs{{\boldsymbol \zeta} f(x-yv)}\; dy
\end{equation}
 See \s.maxfn/ for a discussion of norm bounds for this operator. 
For both this operator, and $\operatorname H^*$, the estimate of weak square integrability is sharp,
as was pointed out to us by M.~Christ \cite{christ}. 
This argument may summarized as follows.  Begin with  a Schwartz function 
$\varphi\ge0$  with frequency support in a small ball about the origin
in the plane.  Then consider $f(x_1,x_2):=\operatorname e^{ix_2}\varphi(x_1,x_2)$.  For any point in the
plane $x=(x_1,x_2)$ with $\abs {x_1}>2$ and $\abs{x_2}<\frac1{100}\abs{x_1}$, 
consider the line that passes through $x $ and the origin.  The real part of $f $ will not change sign on this line, so 
we compute the Hilbert transform with no cancellation.  That is 
 $\operatorname H^* f(x)\simeq{}\operatorname M^*f(x)\simeq\abs{x}^{-1}$.
 And $\abs{x}^{-1}$  is just in weak $L^2$.  See \f.example/.

 \begin{figure}
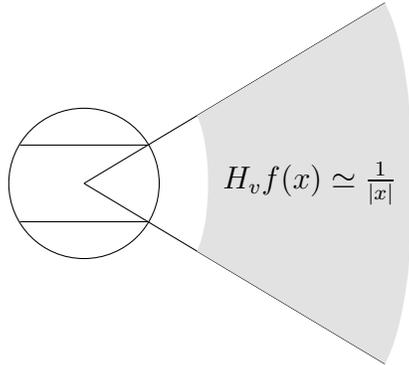
 \label{f.example} 
\begin{center} 
\begin{pgfpicture}{0cm}{0cm}{6cm}{4.5cm}
  	\begin{pgftranslate}{\pgfxy(0,2)}
		\pgfzerocircle{1cm}  \pgfstroke
	{\pgfmoveto{\pgfxy(-.85,.51)}
		\pgflineto{\pgfxy(.85,.51)}   \pgfstroke
		\pgfmoveto{\pgfxy(-.85,-.51)}
		\pgflineto{\pgfxy(.85,-.51)} \pgfstroke}
		 {
		\pgfmoveto{\pgfxy(0,0)}   \pgflineto{\pgfxy(4,2.4)} \pgfstroke 
		\pgfmoveto{\pgfxy(0,0)} \pgflineto{\pgfxy(4,-2.4)} \pgfstroke 
		
		}
		\color{lightgray}
		\begin{colormixin}{45!white}
			\pgfmoveto{\pgfxy(1.5,.9)}  \pgflineto{\pgfxy(4,2.4)} 
			\pgfcurveto{\pgfxy(4.5,1.5)}{\pgfxy(4.5,-1.5)}{\pgfxy(4,-2.4)}
			\pgflineto{\pgfxy(1.5,-.9)}
			\pgfcurveto{\pgfxy(1.7,-.5)}{\pgfxy(1.7,.5)}{\pgfxy(1.5,.9)} 
			\pgffill  
		\end{colormixin}
		\color{black}
		\pgfputat{\pgfxy(3,0)}{\pgfbox[center,center]{$H_v f (x)\simeq\frac1{| x|}  $}}
	\end{pgftranslate}
	\end{pgfpicture}
\end{center} 
 \caption{The function $f $ is depicted by the circle centered at the origin 
 with a planar wave in the vertical direction indicated 
 by the two horizontal lines.  In the shaded region, $H_v f \simeq{}\frac1{\abs x} $.}
\end{figure} 

Whereas, Bourgain's argument for the maximal estimate is not difficult, the Theorem above is of necessity somewhat harder, as it implies the 
pointwise convergence of Fourier series in one dimension. This is an observation in the style of De Leeuw. One considers the trace of the operator in frequency variables along any line in the annulus $1\le\abs \xi\le2$.
 This is Carleson's theorem, \cite{carleson}, but also see \cite{fefferman}. 
As such, we use a method which is adopted from  the proof of Carleson's Theorem given by M.~Lacey and C.~Thiele \cite{laceythiele}.

Perhaps the principal novelty of this paper is the suitable adaptation of the time frequency analysis of Lacey and Thiele to 
the current setting of the plane and measurable choice of directions. 
In comparison to the proof of Lacey and Thiele, we find  a rather precise analogy between the proofs on the real line, and the proofs 
on the plane.  Some differences arise from  the notion of a {\em tile}, which requires some care to define, and the proof 
of the ``size Lemma,'' an orthogonality statement, requires a small amount of innovation.

Finally, there is an outstanding question, attributable to E.~M.~Stein \cite{stein}, concerning the boundedness of the 
Hilbert transform on families of lines that are determined by say a Lipschitz map.  Thus, for a map $v\mid \ZR^2\to 
\{\abs x=1\}$, one wishes to know if 
\begin{equation*}
\int_{-1}^1 f(x-yv(x))\; \frac{dy}y
\end{equation*}
is a bounded operator on say $L^2$.  Positive results for analytic and real analytic vector fields are due to respectively 
Nagel, Stein and Wainger \cite{MR81a:42027}, and to Bourgain \cite{bourgain}.
In a  subsequent paper, the authors \cite{laceyli} will prove that the operator above is bounded on $L^2$ if $v\in C^{1+\ze}$, for any positive \ze.  The results of this paper are a crucial aspect of the proof of this result.

\medskip 

The presentation of this paper has been substantially improved by skilled and generous referees.

\section{Definitions and Principle Lemma}  \label{s.definitions}

\label{s.model}
We begin with some conventions.     We do not keep track of the value of generic absolute constants, instead using the  notation $A\seq{}B$ iff $A\le{}KB$ for some constant $K$.   And $A\simeq B$ iff $A\seq{}B$ and $B\seq{}A$.  We use the notation
$\ind A$ to denote the indicator function of the set $A$.  And the Fourier
transform on $\ZR^2$ is denoted by $\widehat f(\zx)=\int_{\ZR^2} \operatorname e^{-2\zp i x\cdot \zx}f(x)\;dx$, with
a similar definition on the real line.    We use the notation 
\md0
\dashint_A f\;dx:=\abs{A}^{-1}\int_A f\;dx.
\emd

\bthm d.grid  A {\em grid} is a collection of intervals $\z0G$ so that for all $I,J\in\z0G$, we have 
$I\cap J\in\{\emptyset, I, J\}$.  The dyadic intervals are a grid. 
A grid $\z0G$ is {\em central} iff for all $I,J\in\z0G$, with $I\subset_{\not=}J$ we have $100I\subset J$. 
\ethm

The reader can find the details on how to construct such a central 
grid structure in \cite{gl1}.

Let $\zr$ be rotation on $\ZT$ by an angle of $-\zp/2$.  Coordinate axes
for $\ZR^2$ are a pair of unit orthogonal vectors $(e,e_\perp)$  with $e=\zr e_\perp$.  

\bthm d.rectangle
We say that $\zw\subset \ZR^2$ is a {\em rectangle} if it is a product of
intervals with respect to a choice of axes $(e,e_\perp)$ of $\ZR^2$.
We will say that $\zw$ is an {\em annular rectangle}
 if   $\zw=(-2^{l-1},2^{l-1})\times(a,2a)$ for an integer $l$ with
$2^l<a/8$, with respect to the axes $(e,e_\perp)$.
The dimensions of $\zw$ are said to be $2^{l}\times a$.  
 Notice that the face $(-2^{l-1},2^{l-1})\times a $ is tangent to the circle
$\abs \zx=a$ at the midpoint to the face, $(0,a )$.
We say that the {\em scale of}
$\zw$ is $\scl \zw:=2^l$ and that the {\em annular parameter of } $\zw$ is 
$\prm \zw :=a$.
  In referring to the coordinate axes of an annular
rectangle, we shall always mean  $(e,e_\perp)$ as above.  
\ethm

It is imperative to keep in mind that the choice of basis in the definition of a tile depends 
upon the tile in question.   These definitions are illustrated in \f.tile/.  Also see \f.tree/.

 Annular rectangles will decompose our functions in the frequency variables.  But 
our
methods must be sensitive to spatial considerations;  it is this and the
uncertainty principle that motivate the next definition.

\bthm d.dual  Two rectangles $R$ and $\mathsf R$ are said to be {\em dual  }
if they are rectangles with respect to the same basis  $(e,e_\perp)$, thus
$R=r_1\times r_2$ and ${\mathsf R}= {\mathsf r}_1\times {\mathsf  r}_2$ for
intervals $r_i, {\mathsf r}_i$, $ i=1,2$.  Moreover, 
$1\le{}\abs{r_i}\cdot\abs{{\mathsf r}_i}\le2$ for $ i=1,2$.   The product of two dual rectangles we shall refer to 
as a {\em phase rectangle.}  The first coordinate of a phase 
rectangle  we think of as a frequency component and
the second as a spatial component.     
\ethm

\begin{figure}
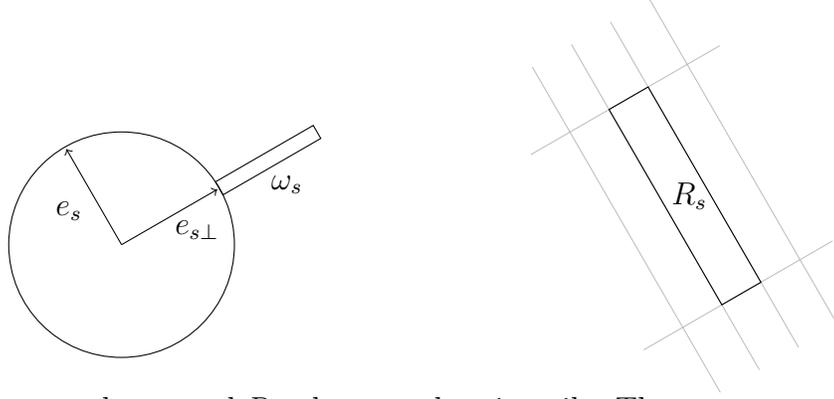
  
  \begin{pgfpicture}{0cm}{0cm}{5cm}{5cm}
%
%
	\begin{pgftranslate}{\pgfpoint{0cm}{1.4cm}}
	\pgfcircle[stroke]{\pgfpoint{0cm}{0cm}}{1.5cm} 
		\begin{pgfrotateby}{\pgfdegree{30}}
		\pgfrect[stroke]{\pgfpoint{1.5cm}{-.1cm}}{\pgfpoint{1.5cm}{0.2cm}}
		
		\pgfsetendarrow{\pgfarrowto}
		\pgfmoveto{\pgfpoint{0cm}{0cm}}\pgflineto{\pgfpoint{1.45cm}{0cm}}\pgfstroke
		\pgfmoveto{\pgfpoint{0cm}{0cm}}\pgflineto{\pgfpoint{0cm}{1.45cm}}\pgfstroke
		
		\end{pgfrotateby}
	\pgfputat{\pgfpoint{2.2cm}{.8cm}}{\pgfbox[center,center]{$\omega_s$}}
	\pgfputat{\pgfpoint{1cm}{.2cm}}{\pgfbox[center,center]{$e_{s\perp}$}}
	\pgfputat{\pgfpoint{-.7cm}{.45cm}}{\pgfbox[center,center]{$e_{s}$}}
	
	\end{pgftranslate}
	\begin{pgftranslate}{\pgfpoint{8.5cm}{.9cm}}
		\begin{pgfrotateby}{\pgfdegree{30}}
			{\color{lightgray}
			\pgfgrid[step={\pgfpoint{0.6cm}{3cm}}]{\pgfpoint{-1.8cm}{-1cm}}{\pgfpoint{1.1cm}{4cm}}
			}		
		\pgfrect[stroke]{\pgfpoint{-.6cm}{0cm}}{\pgfpoint{0.6cm}{3cm}}
		\end{pgfrotateby}
	\pgfputat{\pgfpoint{-.95cm}{1.15cm}}{\pgfbox[center,center]{$R_s$}}
	\end{pgftranslate}
	\end{pgfpicture}
	\caption{ The two rectangles $\omega_s $ and $R_s $ whose product is a tile. The gray rectangles 
	are other possible locations for the rectangle $R_s $.} 
	\end{figure}
   

We consider collections of phase rectangles $\z0{AT}$ which satisfy these
conditions.   For $s,s'\in\z0{AT}$ we write $s=R_s\times \zw_s$, and  require that  
\md5
\Label e.R1 
\text{$\zw_s$ is an annular rectangle }
\\\Label e.R2 
\text{$R_s$ and  $\zw_s$ are  dual   rectangles,}
\\\text{The rectangles $R_s $ are the product of intervals from a central grid.} 
\label{e.R-central}
\\\Label e.R3.5
\{R \mid R\times \zw_s\in\z0{AT}\}\quad\text{partitions $\ZR^2$, for all $\omega_s$.}
\\\Label e.R3 
\prm{\zw_s}={\mathsf {ann}}\quad \text{for some fixed ${\mathsf {ann}}$,}
\\  \Label e.R4 
\sharp\{\omega_s\mid \scl s=2^l,\ \prm s={\mathsf {ann}}\}\ge{}c \operatorname{\mathsf {ann}}2^{-l},\qquad l\in\ZZ.
 \emd
We assume that
there are auxiliary sets $\bw{},\bw{1},\bw2\subset\ZT$ associated to $s$---or
more specifically $\zw_s$---which satisfy these properties.
\md5
\Label e.bw1 
\text{$\z2 \zW:=\{\bw{},\bw{1},\bw2\mid s\in\z0{AT}\}$ is a grid in $\ZT$,}
\\\Label e.bw2 
\bw{1}\cap\bw2=\emptyset,\qquad  \bw{}:=\text{hull}(\bw{1},\bw2),
\\\Label e.bw2/3
\text{$\bw1$ lies clockwise from $\bw2$ on $\ZT$,}
\\\Label e.bw3 
\abs{\bw{}}\le{}K\frac{\scl{\zw_s}}{\prm {\zw_s}},
\\\Label e.bw4 
\{\tfrac\zx{\abs{\zx}}\mid \zx\in\zw_s\}\subset\zr\bw1,
\\\Label e.bw5
\text{to each $\bw{}$ there is a $\z8{\z2 \omega}_s\in \Omega$ such that $2\bw{}\subset \z8{\z2 \omega}_s\subset 4\bw{}$.}
\emd
 Recall that $\zr$ is the rotation that takes $e$ into $e_\perp$. Thus,
 $e_{\zw_s}\in\bw1$.     See \f.omega/.
 



	
	\begin{figure}
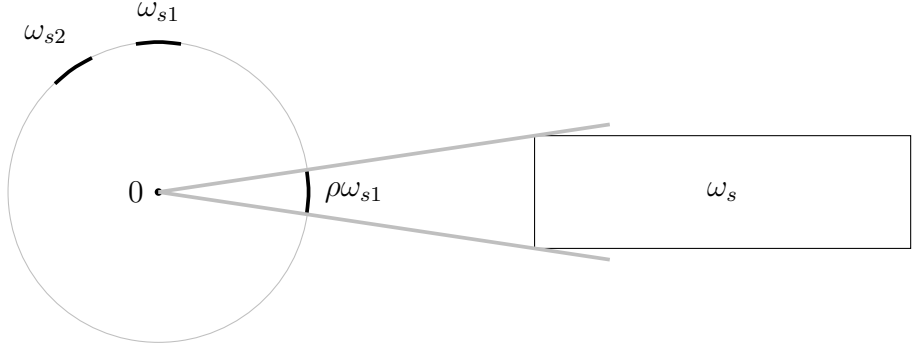
\label{f.omega} 
	
	 \begin{pgfpicture}{0cm}{0cm}{8cm}{4.5cm}
	\begin{pgftranslate}{\pgfxy(-.15,2)}

	{\color{lightgray}  
		\pgfcircle[stroke]{\pgfxy(0,0)}{2cm}
	}
	
	\pgfcircle[fill]{\pgfxy(0,0)}{0.05cm}
	\pgfputat{\pgfxy(-.3,0)}{\pgfbox[center,center]{$0$}}

	\pgfrect[stroke]{\pgfxy(5,-.75)}{\pgfxy(5,1.5)}

		\pgfputat{\pgfxy(2.6,0)}{\pgfbox[center,center]{$\rho \omega_{s1} $}}
	\pgfputat{\pgfxy(0,2.4)}{\pgfbox[center,center]{$\omega_{s1} $}}
	\pgfputat{\pgfxy(-1.5,2.1)}{\pgfbox[center,center]{$\omega_{s2} $}}
	
	\pgfputat{\pgfxy(7.5,0)}{\pgfbox[center,center]{$\omega_s $}}
	
	\pgfsetlinewidth{1.4pt}
	\pgfmoveto{\pgfxy(1.97,.3)}
	\pgfcurveto{\pgfxy(2.005,0)}{\pgfxy(2.005,0)}{\pgfxy(1.97,-.3)}\pgfstroke
	\pgfmoveto{\pgfxy(-.3,1.97)}\pgfcurveto{\pgfxy(0,2.005)}{\pgfxy(0,2.005)}{\pgfxy(.3,1.97)} \pgfstroke
	\begin{pgfrotateby}{\pgfdegree{35}}
	\pgfmoveto{\pgfxy(-.3,1.97)}\pgfcurveto{\pgfxy(0,2.005)}{\pgfxy(0,2.005)}{\pgfxy(.3,1.97)} \pgfstroke
	\end{pgfrotateby}

	{\color{lightgray}  
			\pgfmoveto{\pgfxy(0,0)} 
			\pgflineto{\pgfxy(6,.9)} 	\pgfstroke
			\pgfmoveto{\pgfxy(0,0)} 
			\pgflineto{\pgfxy(6,-.9)}	\pgfstroke
	}

	\end{pgftranslate}
	\end{pgfpicture}
	 \label{f.tile}
	\caption{An annular rectangular $\omega_s $, and three associated subintervals of $\rho \omega_{s1} $, $\omega_{s1} $, 
	and $\omega_{s2} $. }
	\end{figure}
	

Note that 
$\abs{\bw{}}\ge\abs{\bw1}\ge\scl{\zw_s}/\prm{\zw_s}$.  Thus, $e_{\zw_s}$ is in
$\bw{1}$, and  $\bw{}$ serves as ``the angle of uncertainty  associated to
$R_s$.''
{Let us be more precise
about the geometric information encoded into the angle of uncertainty.  Let $R_s=r_s\times
r_{s\perp}$ be as above. Choose another set of coordinate axes $(e',e'_\perp)$ with
$e'\in\bw{}$ and let $R'$ be the product of the intervals $r_s$ and $r_{s\perp}$ in the new coordinate
axes.  Then $K^{-1}_0R'\subset R_s\subset K_0R'$ for an absolute constant $K_0>1$.}

We say that {\em annular tiles} are  collections $\z0{AT}({\mathsf {ann}})$ satisfying the
conditions \e.R1/---\e.bw4/ above.  The constant ${\mathsf {ann}}>0$ appears in  \e.R3/ and \e.R4/. We
extend the definition of $e_\zw$, $e_{\zw\perp}$, $\prm{\zw}$ and $\scl{\zw}$ to
annular tiles in the obvious way, using the notation $e_s$, $e_{s\perp}$, $\prm s $ and
$\scl s $.

\medskip

A phase rectangle will have two distinct functions
associated to it.   In order to define these functions, set 
\md2
\operatorname{Tran}_y f(x)&{}:=f(x-y),\quad y\in\ZR^2\quad \text{(Translation)}
\\
\operatorname{Mod}_\zx f(x)&{}:=\operatorname e^{i\zx\cdot x}f(x),\quad \zx\in\ZR^2\quad \text{(Modulation)}
\\
\operatorname{Dil}_{R_1\times R_2}^p f(x_1,x_2)&{}:= \frac1{(\abs{R_1}\abs{R_2})^{1/p}} 
f\Bigl(\frac{x_1}{\abs{R_1}},\frac{x_2}{\abs{R_2}}\Bigr),\quad \text{(Dilation)}.
\emd
In the last display, $0<p<\infty $,  $R_1\times R_2$ is a rectangle, and the coordinates $(x_1,x_2)$ are those of the rectangle.
Note that the definition depends only on the side lengths of the rectangle, and not the location.  
And that it preserves $L^p$ norm.

For a function $\zvf$ and tile $s\in\z0{AT}$ set 
\md0
 \zvf_s:=\operatorname{Mod}_{c(\zw_s)}\operatorname{ T}_{c(R_s)}\operatorname{Dil}^2_{R_s} \zvf,
\emd
where $c(J)$ is the center of $J$. 
Below, we shall  consider $\zvf$ to be a Schwartz function for which
$\widehat\zvf\ge0$ is supported  in a small ball about the origin in $\ZR^2$, and is identically $1$ on another smaller
ball around the origin.
  
We introduce the tool to decompose the singular integral kernels, and a measurable function $v\mid \ZR^2\to\{\abs x=1\}$ 
that achieves  the maximum in our operator, up to a factor of say $\frac12 $.  We consider a class of functions $\psi_t $, $t>0 $, so that 
\md5
\label{e.zc-Fourier}  \text{Each $\psi_t $ is supported in frequency in $[-\theta-\kappa, -\theta+\kappa ] $. } 
\\
\label{e.zc-Space}  \abs {\psi_t(x) } {}\lesssim{}(1+\abs x )^{-1/\kappa} 
\emd
 In these conditions, $0<\kappa <1$ is a small  fixed  constant that need not concern us.   In the top line, $\theta $ is 
 a fixed positive constant, that is only needed to ensure that \e.dfs2/ below is true.  For the current paper, we will 
need only single choice of function $\psi $, but for the authors'  subsequent paper \cite{laceyli}, countably many 
choices of $\psi $ will be needed.

Define 
\md3\Label e.dfs2 
	\begin{split}
\zf_s(x):={}&\inr \varphi_s(x-yv(x))\zc_s( y)\;dy
\\{}={}&  
\ind {\bw2}(v(x))\inr \varphi_s(x-yv(x))\zc_s( sy)\;dy.
	\end{split} 
\\  \label{e.zc_s} 
\psi_s(y ):={}&\scl s \psi_{\scl s}(\scl s y).
\emd
An essential feature of this definition is that the support of the integral  is contained in the set 
${\{v(x)\in\bw2\}}$, a fact which can be routinely verified, for an appropriate choice of $\theta $ in \e.zc-Fourier/. 
  See \f.localize-v/. That is, we can insert the indicator $\ind {\bw2}(v(x))$ without loss of generality.
The set $\bw2$ serves to localize the vector field, while 
$\bw1$ serves to identify the location of $\zvf_s $ in the frequency
coordinate.

\begin{figure}
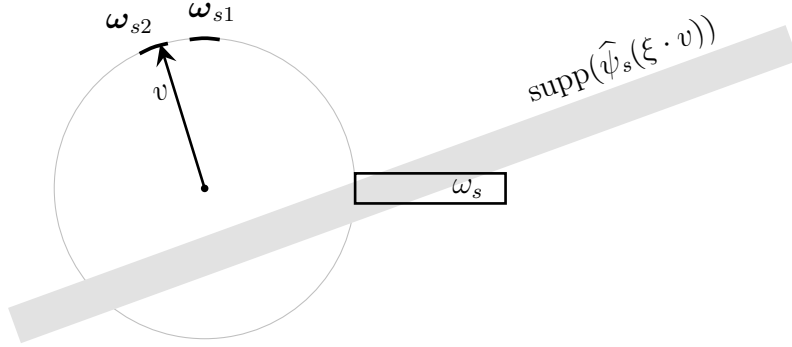
  \label{f.localize-v}
\begin{center} 
		 \begin{pgfpicture}{0cm}{0cm}{8cm}{4.5cm}
	\begin{pgftranslate}{\pgfxy(1.5,2)}

	{\color{lightgray}  
		\pgfcircle[stroke]{\pgfxy(0,0)}{2cm}
	}
	
	\pgfcircle[fill]{\pgfxy(0,0)}{0.05cm}

	\pgfsetlinewidth{1.4pt}

	\pgfmoveto{\pgfxy(-.2,1.98)}\pgfcurveto{\pgfxy(0,2.005)}{\pgfxy(0,2.005)}{\pgfxy(.2,1.98)} \pgfstroke
	\begin{pgfrotateby}{\pgfdegree{20}}
	\pgfmoveto{\pgfxy(-.2,1.98)}\pgfcurveto{\pgfxy(0,2.005)}{\pgfxy(0,2.005)}{\pgfxy(.2,1.98)} \pgfstroke
	\pgfsetendarrow{\pgfarrowsingle} 
	\pgfsetlinewidth{1pt}
	\pgfmoveto{\pgfxy(0,0)}\pgflineto{\pgfxy(.1,1.97)}\pgfstroke
	\pgfputat{\pgflabel{.7}{\pgfxy(0,0)}{\pgfxy(.1,1.97)}{5pt}}{\pgfbox[center,center]{$v$}}
	\color{lightgray} 
	\begin{colormixin}{45!white}
	\pgfrect[fill]{\pgfxy(-3,-1.1)}{\pgfxy(11,.5)}
	\end{colormixin}
	\color{black}	
	\pgfputat{\pgflabel{.8}{\pgfxy(-3,-.6)}{\pgfxy(8,-.6)}{10pt}}{\pgfbox[center,center]{$\text{supp}(\widehat{\psi}_s( \xi\cdot v))$}}

	\end{pgfrotateby}	
	
	\pgfsetlinewidth{1pt}
	\pgfrect[stroke]{\pgfxy(2,-.2)}{\pgfxy(2,.4)}

	\pgfputat{\pgfxy(0.1,2.35)}{\pgfbox[center,center]{$\boldsymbol\omega_{s1} $}}
	\pgfputat{\pgfxy(-1,2.2)}{\pgfbox[center,center]{$\boldsymbol\omega_{s2} $}}
	
	\pgfputat{\pgfxy(3.5,0)}{\pgfbox[center,center]{$\omega_s $}}

	\end{pgftranslate}
		\end{pgfpicture}
	\end{center} 
	
	\caption{ The support of $\phi_s $ is contained in $v^{-1}(\boldsymbol\omega_{s2}) $. } 
	\end{figure}

\smallskip 

{\bf Remark:}  It seems likely that we could prove our theorem using the simpler definition $\zf_s=(\ind {\bw2}\circ v )\zvf_s$. 
This would be the natural analog of the decomposition used by Lacey and Thiele \cite{laceythiele}. 
But, in our subsequent paper on this subject \cite{laceyli}, we will need to consider a truncation of the Hilbert transform, which suggests the definition of $\zf_s$ we have adopted above.  We will also want to 
rely upon facts proved in this paper, to deduce our theorems concerning smooth vector fields. 

\smallskip

The model operators we consider are  defined by 
\md0
\z0C_{\mathsf{ann},v}f:=  \sum_{\substack{s\in\z0{AT}({\mathsf {ann}})} }\ipf \phi_s\,   .
 \emd
 In this display, $\z0{AT}({\mathsf {ann}}):=\{s\in\z0{AT}\mid \prm s={\mathsf {ann}}\}$.


 \bthm  l.model   Assume that the vector field is only measurable.  The operator $\z0C_{\mathsf{ann},v}$ extends to a 
 bounded map  from $L^2$ into weak $L^2$, 
 and $L^p$ into itself for $2<p<\infty$.  The norm of the operator is independent of ${\mathsf {ann}}$.
\ethm

We shall  prove that $\z0C_{\mathsf{ann},v}$ maps $L^2$ into weak $L^2$, with constant independent of $\mathsf{ann}$
and $v$.  By duality, it suffices to show that for all $f\in L^2\cap L^\infty$ of $L^2$ norm one, and sets $F\subset\ZR^2$ 
of finite measure
\begin{equation}\label{e.cjabs}
\abs{ \ip \z0C_{\mathsf{ann},v} f,\ind F.}\le{}\sum_{s\in\z0{AT}({\mathsf {ann}})}\abs{\ipf \ip \phi_s,\ind F.}\seq\abs{F}^{1/2}.
\end{equation}
By dilation invariance of the $\z0C_{\mathsf{ann},v}$ with respect to powers of $2$, we can further take $1\le{}\abs F<2$.

For the case of $L^p$, $2<p<\zI$, we shall demonstrate the restricted type estimate
\md0
\abs{\{\z0C_{\mathsf{ann},v} \ind E >\zl\}}\lesssim\zl^{-p}\abs{E},\qquad 2<p<\zI.
\emd
We need only consider this for $\zl>1$, by the weak $L^2$ bound.  Moreover, this inequality follows from 
\md3 \nonumber
\abs{\ip \z0C_{\mathsf{ann},v} \ind E,\ind F.}\le{}&\sum_{s\in\z0{AT}({\mathsf {ann}})}\abs{\ip \zvf_s, \ind E. 
\ip \phi_s,\ind F.}
\\{}\lesssim{} &
\abs{F}\{\log (\abs E/\abs F)\},\qquad 9\abs F<\abs E.\label{e.p>2}
\emd
Observe that by dilation invariance, we can assume that $\abs E\simeq1$.  We also assume that the vector field $v$ is defined
only on $F$. 

The proofs of \e.cjabs/ and \e.p>2/  are given in \s.outline/, and follow lines of argument that are similar to the one dimensional 
case, as given by Lacey and Thiele \cite{laceythiele} for the case $p=2 $, and as in the work of Grafakos, Tao, and Terwilleger 
\cite{gtt} 
for the case of $p>2 $.  Also see the survey by Lacey \cite{esi}.    

 {\bf Remark.} The results for $p>2 $ are available essentially because functions that are in $L^p $ are also locally square integrable.  In one 
dimension, the proof of convergence of Fourier series for $1<p<2 $ is closely linked to the boundedness of the Hardy Littlewood 
maximal function.  In the current setting, the corresponding maximal function would be the unbounded Kakeya maximal function. 
The relationship between maximal functions and $L^p $ bounds for $1<p<2 $ has a subtle role to play in our 
companion paper \cite{laceyli}.

\section{A Maximal Function Estimate } \label{s.maxfn} 

We outline a proof of the weak type estimate on $L^2 $ for the directional maximal function $\operatorname M^* $ given in 
\e.dir-max/.   

\bthm p.maxfn   For functions $f\in L^2 $, and $\lambda>0 $, we have the estimate 
\md0
\abs{ \{ \operatorname M^* f>\lambda\} } {}\lesssim{}\lambda^{-2}{ \norm f.2.^2}
\emd
\ethm

It suffices to give the same bound for the model operator
\md0
\z0M_{\mathsf{ann},v}f:=  \sup_{\substack{s\in\z0{AT}({\mathsf {ann}})} }\abs {\ipf \phi_s }   .
\emd
It is a remark that this inequality is an immediate consequence of \l.model/.  
The sum defining $\z0C _{\mathsf{ann},v} $ is in fact unconditional  convergent as a sum over tiles. 
As a consequence the square function below, which trivially dominates $\z0M_{\mathsf{ann},v} $, satisfies the weak type inequality.  
\md0
\Bigl[ \sum_{s\in\z0{AT}({\mathsf {ann}})} \abs {\ipf \phi_s }  ^2\Bigr]^{1/2}
\emd

There is however a more direct proof.  The fundamental fact, a simplier instance  of the Size Lemma, \l.size/ above, 
is this. 

\bthm l.size-simple   Fix $f\in L^2 $ and $\lambda>0 $.  Let $\z0S\subset \z0{AT}({\mathsf {ann}}) $ be a collection of tiles which are pairwise 
incomparable with respect to the partial order `$<$,' and moreover for  each $s\in\z0S $, $\abs{\ipf}\ge \lambda\sqrt{\abs{R_s}} $. 
Then it is the case that 
\md0
\sum _{s\in\z0S} \abs{ R_s }\le{} \lambda ^{-2 }\norm f.2.^2 
\emd
\ethm 

This is proved in e.g.~Barrionuevo and Lacey \cite{MR1955263}, or one can inspect the proof of \l.size/ for a number of 
simplifications which apply in this case.  

It is then clear that the maximal function below maps $L^2 $ into weak $L^2 $. 
\md0
\operatorname  A_1 f :=\sup _{\z0{AT}({\mathsf {ann}}) }\frac {\abs{\ipf } }{ \sqrt{\abs{R_s}} }\mathbf 1 _{R_s}
\emd
For $t>2 $, set 
\md0
\operatorname  A_t f :=\sup _{\z0{AT}({\mathsf {ann}}) }\frac {\abs{\ipf } }{ \sqrt{\abs{R_s}}} \mathbf 1 _{tR_s}
\emd
This also maps  $L^2 $ into weak $L^2 $ with norm at most $ {}\lesssim{}t^2 $.  
Indeed, given $\lambda>0 $, let $\z0S $ denote the maximal tiles with $\abs{\ipf}\ge \lambda\sqrt{\abs{R_s}} $.
Then, it is the case that 
\md2
\{ \operatorname A_t f> \lambda \}\subset\bigcup _{s\in\z0S} KtR_s
\emd
for an absolute constant $K $.  Hence, the claimed bound follows from \l.size-simple/. 
Then observe that 
\md0
\z0M_{\mathsf{ann},v}f
{}\lesssim{}\sum _{k=0}^\infty 2^{-4k} \operatorname A _{2^k} f . 
\emd

 \section{Proof of \t.main/  }

 We show that \l.model/ implies our main result, \t.main/.  Indeed, the Hilbert transform in the direction of $v $ is, in the limit, 
 an appropriate average of the discrete sums formed in the previous section.  Let 
 $\zl$ be a smooth radial function satisfying 
\md1 \label{e.zk}
 \ind {[\frac54,\frac74]}(\abs \zx)\le\widehat\zl\le\ind {[1,2]}(\abs \zx)
\emd
Let $\zl_t(y)=t^{2}\zl(ty)$, where we choose the dilation so that $t $ is the Fourier parameter, as opposed to the more 
common definition $t^{-2}\zl(y/t) $ where $t $ is the spatial parameter.

Let $K$ be the distribution on \ZR 
\md0
\sum_{j=-\zI}^\zI 2^j \zc(2^j y)
\emd
 where $\psi $ is a Schwartz function with  $\widehat\zc \ge0$  supported in a small neighborhood of $-1$.  In particular, 
the distribution 
\md0
\int_0^1 2^s K(2^s y)\; ds
\emd
is a non--zero multiple of the distribution associated with projection onto the negative frequencies on \ZR, which is itself
 a linear combination of the identity operator and the Hilbert transform. 

To prove \t.main/ it suffices to demonstrate the same norm estimates for the linear operators  in which the 
distribution $1/y $ is replaced by $K(y) $, namely for 
\md1 \label{e.Tv}
\operatorname T_v f(x)=\int \zl_{\mathsf {ann}}*f (x-y v(x))K(y)\; dy 
\emd
in which ${\mathsf {ann}}>0$, and the measurable $v(x)$ are arbitrary.  
 
 \medskip 
 
 For values of $2^j<{\mathsf {ann}}$, let 
 \md2
 \z3{\operatorname S}_{j,{\mathsf {ann}}}f=\sum_{\substack{s\in\z0{AT}({\mathsf {ann}}) \\ \scl s=2^j}} \ipf \zvf_s 
 \\
 \operatorname T_{j,v}f(x)=\int f(x-yv(x))2^j\zc(2^j y)\; dy 
 \emd
 Also, let $\operatorname{Rot}_\zt$ be the operation of rotation by angle $\zt$. The main point 
concerns the operator 
\md1 \label{e.avg}
\begin{split} 
\operatorname S_{j,{\mathsf {ann}}}:=\lim_{Y\to\zI}\dashint_{\text{Box}(Y)}  & \operatorname {Dil}_{[0,2^{-s})\times[0,2^{-s})}^2 
\operatorname{Rot}_{-\zt}\operatorname{Tran}_{-y}
\\  
&  \z3{\operatorname{S}}_{j,{\mathsf {ann}}} 
\operatorname{Tran}_{y}\operatorname{Rot}_{\zt} \operatorname {Dil}_{[0,2^s)\times[0,2^s)}^2 \;dy\,d\zt\, ds
\end{split} 
\emd
 where $\text{Box}(Y)=[0,Y]^2\times[0,2\zp]\times[-1,1]$.  Observe these points.  (1) The integral is an average over dilations, rotations, and translations. 
 Note that the averaging over dilations is done uniformly in all directions, and that we are implicitly averaging with 
 respect to the Haar measure on $\mathbb R_+ $, with the rationale being that the dilations we are using are a representation 
 of that group.  (2) Due to the assumption \e.R3.5/,  the limit in $Y $ will exist.  If $\omega_s $ were fixed
we would only need to average over all translations indexed by a fixed rectangle.  (3) The rectangles $\omega_s $ fill 
a fixed percentage of the annulus, due the assumption \e.R4/.    (4) The  limit, applied to 
 a Schwartz function $f$ is seen to exist.

 \bthm l.2convolve 
 For each $2^j<{\mathsf {ann}}$, we have the identity 
 \md0
\operatorname  S_{j,{\mathsf {ann}}}\zl_{{\mathsf {ann}}}*f=c(j,{\mathsf {ann}})\zl_{{\mathsf {ann}}}*f
 \emd
 where the constant $c(j,{\mathsf {ann}})$satisfies $c^{-1}<\abs{c(j,{\mathsf {ann}})}<c$, for some absolute constant $c$. 
 \ethm 
 
 To deduce bounds for \e.Tv/, observe that \l.model/ concerns norm bounds for the sums 
 \md0
 \sum_j c(j,{\mathsf {ann}})^{-1} \operatorname T_{j,v} \z3{\operatorname{S}}_{j,{\mathsf {ann}}} \zl_{{\mathsf {ann}}} *f 
 \emd 
 Of course the coefficients $c(j,{\mathsf {ann}})^{-1}$ do not appear in \l.model/.  Yet the placement of the 
 absolute values in \e.cjabs/  and \e.p>2/ demonstrates that the sum is unconditional over tiles, so that we can impose 
 an arbitrary bounded sequence of coefficients, as we have done here. 
 
 The same norm bounds hold for averages of these sums, such as the averages that occur in \e.avg/.  Hence, by the
 Lemma just stated, 
 we deduce the norm bounds for the operators in \e.Tv/.  It remains to prove our Lemma.

Before proving \l.2convolve/, we  record a simple lemma on convolutions.  

\bthm l.convolve  Let $\zvf$ and $\zf$ be real valued Schwartz functions on $\ZR^2$.  Then, 
\md2
\int_{[0,1]^2}\sum_{m\in\ZZ^2} 
\ip f, \operatorname {Tran}_{y+m}\zvf. \operatorname {Tran}_{y+m}\zf \; dy=f*\zF
\\
\text{where}\qquad \zF(x)=\iint \overline{\zvf(u)}\zf(x+u)\; du.
\emd
In particular, $\widehat{\zF}=\overline{\widehat\zvf}\widehat \zf$. 
\ethm

The proof is immediate.  The sum and integral in question is in fact the integral 
\md0
\int_{\ZR^2} f(z)\overline{\zvf(z-y)}\zf(x-y)\; dydz
\emd
and one changes variables, $u=z-y$.

\begin{proof}[Proof of \l.2convolve/]

 Fix $\zw\in\zW$, and let 
 \md0
 \z3{\operatorname{S}}_\zw f:=\sum_{\substack{s\in\z0{AT}({\mathsf {ann}})\\ \zw_s=\zw}} \ipf \zvf_s 
 \emd
 Recall \e.R3.5/.  Consider the limit 
 \md0
 \operatorname S_\zw :=\lim_{Y\to\zI}Y^{-2} \dashint_{[0,Y]^2} \operatorname{Tran}_{-y} \z3{
 \operatorname{S}}_\zw \operatorname{Tran}_{y} \; dy
 \emd
 This limit can be explicitly written as an integral which is as in \l.convolve/ with an additional dilation. 
 Thus, this limit 
 is convolution with respect to $\zF_\zw$, where $\widehat\zF_\zw=\widehat \zvf_\zw \overline{\widehat\zvf_\zw}$.  
 
 By choice of $\zvf$, it has non--zero Fourier transform, and is identically one on a small ball around the origin. 
 Hence $\zF_\zw$ has Fourier transform with the same properties.  
 The definition of $\operatorname  S_{j,{\mathsf {ann}}} $ also incorporates an average over dilations and rotations. 
 Concerning the average over rotations, 
 recall \e.R4/, which states that we have an essentially 
 maximal number of $\zw_s$ of a given scale.  
 Concerning the dilations, observe that 
 \md0
 \operatorname {Dil}^2_{[0,t)^2}\bigl[ F *[ \operatorname {Dil}^2_{[0,1/t)^2} G]={} 
 [\operatorname {Dil}^1_{[0,t)^2} F]*G.
 \emd

 These observations lead to the conclusion that 
 \md2
\operatorname  S_{j,{\mathsf {ann}}}f&{}=c(j,{\mathsf {ann}}) \zF*f
\\
\zF&{}=\int_{-1}^1\int_{[0,2\zp]}  \operatorname {Dil}^1_{[0,2^t)^2}\operatorname{Rot}_{\zt}
\zF\; d\zt\,dt
\emd
The leading constant in the top line satisfy the claimed bound  $c<\abs{ c(j,{\mathsf {ann}}) }<c^{-1} $ for 
a positive constant $c $.   The main point here is assumption 
 \e.R4/.  This function \zF\ has Fourier transform that is identically constant 
for all ${\mathsf {ann}}<\abs\zx<2{\mathsf {ann}} $.  The value of this constant depends only on the 
choice of $\zvf $.  
 \end{proof}


\section{The Main Lemmas} \label{s.outline}

We need these notions  associated to tiles and sets of tiles, all directly inspired by the proof of Carleson's theorem  as 
given by 
Lacey and Thiele \cite{laceythiele}.  There are however some changes in the current context, due to the setting of the plane, 
and our frequency analysis in an annulus.  

An important heuristic is that given a tile $s $, the interval ${\z2 \omega_s}$ serves as an ``angle of uncertainty''{} 
for the rectangle $R_s $.   The geometric consequence of this can be described as follows.   Write 
$R_s=R_{s1}\times R_{s2} $ in the coordinate axes $(e_s,e_{s\perp}) $.   There is an absolute constant 
$K $ so that the following holds.  Suppose that $(e,e_\perp) $ is some other 
choice of coordinate axes with 
\md0
\abs{e_s-e},\ \abs{e_{s\perp}-e_\perp}\le{}\abs{\z2 \omega_s}.
\emd
Write $\widetilde R_s =R_{s1}\times R_{s2} $ but the choice of coordinate axes is $(e,e_\perp) $.  Then, we have 
$R_s\subset K \widetilde R_s  $ and $\widetilde R_s \subset K R_s $.   The implication is that $R_s $ and $\widetilde R_s $
are essentially the same rectangles. See \f.uncertain/.

 	\begin{figure}
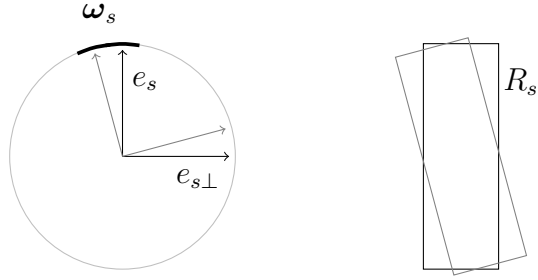
\label{f.uncertain} 
	\begin{center}
	 \begin{pgfpicture}{0cm}{0cm}{5.5cm}{4.5cm}
	\begin{pgftranslate}{\pgfxy(-.15,2)}

	
	\pgfrect[stroke]{\pgfxy(4,-1.5)}{\pgfxy(1,3)}
	\pgfputat{\pgfxy(5.3,1)}{\pgfbox[center,center]{$R_s$}}

	\begin{pgftranslate}{\pgfxy(4.5,0)}
	\color{gray}
	\begin{pgfrotateby}{\pgfdegree{15}}
	\pgfrect[stroke]{\pgfxy(-.5,-1.5)}{\pgfxy(1,3)}
	\end{pgfrotateby}
	\end{pgftranslate}

	{\color{lightgray}  
		\pgfcircle[stroke]{\pgfxy(0,0)}{1.5cm}
	}
	\pgfsetendarrow{\pgfarrowto}	
	\pgfmoveto{\pgfxy(0,0)} \pgflineto{\pgfxy(0,1.4)} \pgfstroke	
			\pgfputat{\pgfxy(.3,1)}{\pgfbox[center,center]{$e_s$}}
	\pgfmoveto{\pgfxy(0,0)} \pgflineto{\pgfxy(1.4,0)} \pgfstroke	
			\pgfputat{\pgfxy(1,-.3)}{\pgfbox[center,center]{$e_{s\perp}$}}
		\begin{pgfrotateby}{\pgfdegree{15}}	
			\color{gray}
			\pgfmoveto{\pgfxy(0,0)} \pgflineto{\pgfxy(0,1.4)} \pgfstroke	
			\pgfmoveto{\pgfxy(0,0)} \pgflineto{\pgfxy(1.4,0)} \pgfstroke	
		\end{pgfrotateby}	
	\pgfclearendarrow	
	
	\pgfputat{\pgfxy(-.3,1.9)}{\pgfbox[center,center]{$\boldsymbol\omega_{s} $}}
	\pgfsetlinewidth{1.4pt}
	\pgfmoveto{\pgfxy(1.477,.225)}
	\pgfmoveto{\pgfxy(-.225,1.477)}\pgfcurveto{\pgfxy(0,1.50375)}{\pgfxy(0,1.50375)}{\pgfxy(.225,1.477)} \pgfstroke
		\begin{pgfrotateby}{\pgfdegree{15}}
		\pgfmoveto{\pgfxy(1.477,.225)}
		\pgfmoveto{\pgfxy(-.225,1.477)}\pgfcurveto{\pgfxy(0,1.50375)}{\pgfxy(0,1.50375)}{\pgfxy(.225,1.477)} \pgfstroke
		\end{pgfrotateby}

	\end{pgftranslate}
	\end{pgfpicture}
	\end{center}
	\caption{On the left, the interval $\boldsymbol\omega_{s} $ and the axes $(e_s,e_{s\perp}) $.  On the right 
	a rectangle $R_s $.  In gray, the axes and rectangle $R_s $ are rotated by an amount less than the length of 
	$\boldsymbol\omega_{s}  $. }
	\end{figure}

	There is a natural partial order on tiles, 
but it does take some care to define. 
For a fixed tile $s$, let $(e_s, e_{s\perp})$ be the coordinate system
associated with $R_s$.   While we have not done so to date, at this point we insist 
that the origin in each such coordinate system specify the same point in the plane. 
Let ${\z0 P} $ be the set of the projections 
of all $R_s$'s onto $e_s$ and ${\z0 P}_\perp $ be the set of the projections 
of all $R_s$'s onto $e_{s\perp}$.  We  have assumed that $\z0P $ is a central 
grid in \e.R-central/.   Namely, recalling \d.grid/,  any two 
intervals $I\not=J$ in ${\z0 P}$ which intersect,   satisfy either $100I\subset J$, or $100J\subset
I$.  The reader can find the details on how to construct such a central 
grid structure in \cite{gl1}.   As all tiles lie in a single annulus, the collection of intervals 
${\z0 P}_\perp $ consist of disjoint intervals of equal length that partition $\mathbb R $.

Let $\pi(R, e_s)$ be the minimal interval 
in  ${\z0 P}$ containing the projection of $R$ onto $e_s$, and 
  let $\pi(R, e_{s\perp})$ be the minimal interval 
in  ${\z0 P}_\perp$ containing the projection of $R$ onto $e_{s\perp}$.

 For two tiles $s, s' $, write  $s<s'$ if and only if  
${\z2 \omega}_s\supset{\z2 \omega}_{s'}$, 
$\pi(R_s,e_{s'})\subset \pi(R_{s'},e_{s'})$, 
and $\pi(R_{s'},e_{s'\perp})=\pi(R_{s'},e_{s'\perp})$.  This is illustrated in \f.partial/.

We  have been careful to define this partial order in such a way that it is transitive and that we have the conclusion
\begin{equation}\label{e.<property}
\text{If $\bw{}\times R_{s}\cap \bw'\times R_{s'}\not=\emptyset$, then  $s$ and $s'$ are comparable under `$<$'.}
\end{equation}
 
\begin{proof}
We can assume that $\bw{}\supset \bw' $, as these intervals of the circle come from a grid.  
And so it suffices to check that $\pi(R_s,e_{s'})\subset \pi(R_{s'},e_{s'})$.  By the assumption \e.R-central/, 
in which $\kappa $ is a small positive constant, we would have $\pi(R_s,e_{s'})\subset 100 \pi(R_s,e_{s}) $, 
and the last set is contained  in $\pi(R_{s'},e_{s'}) $, by the central grid property. 
\end{proof}

\begin{figure}
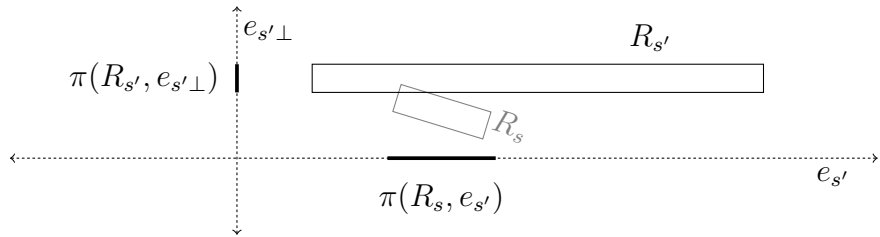
\label{f.partial} 
	\begin{center}
	 \begin{pgfpicture}{-1cm}{0cm}{6.5cm}{4cm}
	\begin{pgftranslate}{\pgfxy(0.3,2)}     
	\begin{pgfrotateby}{\pgfdegree{-17}}
	\begin{pgfscope}
	\color{gray}
	\pgfrect[stroke]{\pgfxy(1.8,.5)}{\pgfxy(1.25,.375)}
	\pgfputat{\pgflabel{.7}{\pgfxy(3.35,.5)}{\pgfxy(3.35,.875)}{0pt}}{\pgfbox[center,center]{$R_{s}$}}
	\end{pgfscope}
	\begin{pgftranslate}{\pgfxy(0,-.7)}
	\begin{pgfrotateby}{\pgfdegree{17}}
%
	\color{black}
	\pgfrect[stroke]{\pgfxy(1,.875)}{\pgfxy(6,.375)}
	\pgfputat 
			{\pgflabel{.9}{\pgfxy(1,1.25)}{\pgfxy(6,1.25)}{10pt}}
					{\pgfbox[center,center]{$R_{s'}$}}
%
%
	\begin{pgfscope} 
	\pgfsetdash{{1pt}{1pt}}{0pt} 
	\pgfsetendarrow{\pgfarrowto}\pgfsetstartarrow{\pgfarrowto}
	\pgfmoveto{\pgfxy(0,-1)}\pgflineto{\pgfxy(0,2)}\pgfstroke
	\pgfputat 
		{\pgflabel{.9}{\pgfxy(0,-1)}{\pgfxy(0,2)}{-12pt}}
			{\pgfbox[center,center]{$e_{s'\perp}$}}
	\pgfputat
		{\pgflabel{.95}{\pgfxy(-3,0)}{\pgfxy(8.5,0)}{-7pt}}
			{\pgfbox[center,center]{$e_{s'}$}}
	\pgfmoveto{\pgfxy(-3,0)}\pgflineto{\pgfxy(8.5,0)}\pgfstroke
			\pgfclearendarrow\pgfclearstartarrow
		\end{pgfscope}
%
		\pgfsetlinewidth{1.4pt}
			\pgfmoveto{\pgfxy(2,0)}\pgflineto{\pgfxy(3.44,0)} \pgfstroke
			\pgfputat
				{\pgflabel{.5}{\pgfxy(2,0)}{\pgfxy(3.44,0)}{-15pt}}
					{\pgfbox[center,center]{$\pi(R_{s},e_{s'})$}}
			\pgfmoveto{\pgfxy(0,.875)}\pgflineto{\pgfxy(0,1.25)}\pgfstroke
			\pgfputat
				{\pgflabel{.5}{\pgfxy(0,.875)}{\pgfxy(0,1.25)}{35pt}}
					{\pgfbox[center,center]{$\pi(R_{s'},e_{s'\perp})$}}
	\end{pgfrotateby}
	\end{pgftranslate}
	\end{pgfrotateby}
	\end{pgftranslate}
	\end{pgfpicture}
	\end{center}
	\caption{ For tiles $s<s' $, only showing possible relative positions of $R_s $ and $R_s' $, as well as the projections 
	that define the partial order. }
	\end{figure}

A {\em tree} is   a collection of tiles $\T\subset\z0{AT}({\mathsf {ann}})$, for which there is a (non--unique) tile
$\zw_\T\times R_\T\in\z0{AT}({\mathsf {ann}})$ with $s<\zw_\T\times R_\T$ for all $s\in\T$.  For $j=1,2$, 
call $\T$ a $j$--tree if the tiles $\{\bw{i}\times R_s \mid s\in\T\}$ are pairwise disjoint. 
$1 $--trees are especially important for us;  \f.tree/ presents an useful illustration of such a tree.

Observe that for any tree $\T $ we have, due to our comments about intervals of uncertainty,  
\md0
R_s\subset K R_\T ,\qquad s\in \T . 
\emd
Thus, the spatial locations of the tiles are localized.  
In addition, as ${\z2 \omega}_s \subset \z2\zw_\T $, we can regard all rectangles $R_s $ as having a \emph{fixed }
set of coordinate axes, namely those of the top.   This is a key point of the ``Tree Lemma''{} below.

Fix a positive rapidly decreasing function $\chi$, and for rectangle $R$ set 
 \md5 \label{e.zqrs}
\chi_{R}=\operatorname {Tran}_{c(R)}\operatorname {Dil}^1_R \chi
\\ \nonumber
\dense s:=\sup_{s<s'} \int_{v^{-1}(\bw')} \zq_{R_{s'}}\; dx,
\\\nonumber
\dense \S:=\sup_{s\in\S}\dense s,
\\\nonumber
\sh \S:=\bigcup_{s\in\S}R_{s} \qquad \text{(the {\em shadow} of \S)}
\\\nonumber
\zD(\T)^2:=\sum_{s\in\T}\frac{\abs{\ipf}^2}{\abs{R_{s}}}\ind {R_{s}},\qquad \text{$\T$ is a $1$--tree,} 
\\\nonumber
\size \S:=\sup_{\substack{ \T\subset \S \\ \text{$\T$ is a $1$--tree}}} \bigl[\abs{\sh \T}^{ -1} 
\sum_{s\in\T} 
\abs{\ipf}^2\bigr]^{1/2}
\emd
  Recall that  we are to prove \e.cjabs/ and \e.p>2/.  In \e.cjabs/,  $\norm f.2.\simeq1$, $\abs F\le{}2$, and
  $v$ is defined only on $F$.

 The first two Lemmas concern density and size, and are the same form.   The third Lemma, the `Tree Lemma' relates 
 size, density and trees. 

\bthm l.dense  
Any $\S\subset\z0{AT}({\mathsf {ann}})$ is the union of $\S_{\text{\rm heavy}}$ and $\S_{\text{\rm light}}$ 
satisfying these conditions.  
\md0
\dense {\S_{\text{\rm light}}}\le{}\tfrac12\dense \S.
\emd
The collection $\S_{\text{\rm heavy}}$ is a union of trees $\T\in\vZT_{\text{\rm heavy}}$, with 
\md1
\label{e.dense}
\sum_{\T\in\vZT_{\text{\rm heavy}}}\abs{\sh\T}\seq{} {\dense \S}^{-1}\abs F.
\emd
\ethm


\bthm l.size  
Any $\S\subset\z0{AT}({\mathsf {ann}})$ is the union of $\S_{\text{\rm big}}$ and $\S_{\text{\rm small}}$
satisfying these conditions.  
\md0
\size {\S_{\text{\rm small}}}\le{}\tfrac12\size \S.
\emd
The collection $\S_{\text{\rm big}}$ is a union of trees $\T\in\vZT_{\text{\rm big}}$, with 
\md1
\label{e.size}
\sum_{\T\in\vZT_{\text{\rm big}}}\abs{\sh\T}\seq{} {\size \S}^{-2}.
\emd
\ethm

\bthm l.tree   
For any tree $\T$ we have the estimate 
\md0
\sum_{s\in\T}\abs{\ipf \ip \zf_s,\ind E.}\seq{}\dense \T \size\T \abs{R_\T}.
\emd
\ethm

The first lemma has a proof which is essentially identical to the proof of the ``mass lemma'' in M.~Lacey and C.~Thiele
\cite{laceythiele}.  We do not give the proof. 
The second lemma follows the well established lines of proof, yet must introduce a new element or two 
to address the two dimensional setting.  The complete proof is given.  The proof of the last lemma is quite close to that 
of M.~Lacey and C.~Thiele  \cite{laceythiele}.  We shall give a proof.

\medskip 

The lemmas are combined in this way to prove \e.cjabs/, and hence \l.model/ in the case of the weak $L^2$ estimate.
   \l.dense/ and \l.size/ should be applied so that their principal estimates 
\e.dense/ and \e.size/ are approximately equal.
 The density of $\z0{AT}({\mathsf {ann}})$ is at most a constant.  The size of $\z0{AT}({\mathsf {ann}})$ 
 is at most a constant times the 
$L^\infty$ norm of $f$.  Thus, we can take the set of all tiles $\z0{AT}({\mathsf {ann}})$ and decompose
it into sub collections $\S_ \sigma $, 
$\sigma\in\{2^n\mid n\in\ZZ\}$,
so that $\S_ \sigma$ is the union of trees $\T\in\vZT_ \sigma$ such that 
\md4
\sum_{\T\in\vZT_ \sigma}\abs{R_\T}\seq{}\sigma,
\\
\dense {\S_ \sigma}\seq{}\min( 1, \sigma^{-1}),
\\
\size {\S_ \sigma}\seq{} \sigma^{-1/2}.
\emd
Hence, it follows that 
\md2
\sum_{s\in\S_ \sigma}\abs{\ipf \ip \phi_s,\ind E.}\seq{}&\sum_{\T\in\vZT}\min( 1, \sigma^{-1})\sigma^{-1/2}\abs{R_\T}
\\
{}\seq{}&\min(\sigma^{1/2},\sigma^{-1/2}).
\emd
This estimate is summable over $\sigma\in\{2^n\mid n\in \ZZ\}$ to an absolute constant. This completes the proof of \e.cjabs/.

 \smallskip
 
 A small variation on this argument proves \e.p>2/.  In this instance, the function $f=\ind E$, with $\abs E\simeq1$, 
 so that the size of $\z0{AT}({\mathsf {ann}})$ is  ${}\lesssim1$.   And $\z0{AT}({\mathsf {ann}})$ is a union of 
sub collections $\S_ \sigma $, 
$\sigma\in\{2^n\mid n\in\ZN\}$, 
so that $\S_ \sigma$ is the union of trees $\T\in\vZT_ \sigma$ satisfying 
\md4
\sum_{\T\in\vZT_ \sigma}\abs{R_\T}\seq{}\sigma\abs F,
\\
\dense {\S_ \sigma}\seq{} \sigma^{-1},
\\
\size {\S_ \sigma}\seq{} \min( 1,(\sigma\abs{F})^{-1/2}).
\emd
Hence, it follows that 
\md0
\sum_{s\in\S_ \sigma}\abs{\ip \mathbf 1 _{E},\zvf_s. \ip \phi_s,\ind F.}\seq{}\abs{F}^{1/2}\min(\abs{F}^{1/2},\sigma^{-1/2}).
\emd
Recall that in this instance, $\abs F\le\frac13$. 
This estimate is summable over $\sigma\in\{2^n\mid n\in \ZN\}$ to $\abs F\abs{\log \abs F}$. This completes the proof of \e.p>2/.
Our proof of \l.model/ is complete, aside from the proofs of the three Lemmas of this section.


\section{The Size Lemma: Orthogonality}

We give the proof of \l.size/.
 We find that the proof involves for the most part a standard argument in the literature.  See for instance the 
proof of the energy lemma in M.~Lacey and C.~Thiele \cite{laceythiele}.
Yet there is a point at which we will rely upon the strong maximal function, with respect to a choice of axes that is
specified in a particular way by the set of tiles.

We can assume that all $\bw{} $ lie in a fixed half circle.  Note that in \tree1s, that as the scale of tiles 
increase, the intervals $\bw1 $ move off in a clockwise direction, as depicted in \f.tree/.  We take advantage of this 
in a specific way below.

The initial step of the proof is to construct a collection of  $1$--trees $\T\in\vZT_{+}$ and use them to construct the 
collection of trees $\vZT_{\text{\rm big}}$.  The process is inductive. Let $\size \S=\sigma$.  Initialize 
\begin{equation*}
\vZT_{+}:=\emptyset,\qquad \vZT_{\text{\rm big}}=\emptyset,\qquad \S^{\text{\rm stock}}=\S.
\end{equation*}
While  $\size {\S^{\text{\rm stock}}}>\sigma/2$,   select 
a $1$--tree $\T\in\S^{\text{\rm stock}}$ such that 
\begin{equation}\label{e.treeselect}
\sum_{s\in\T}\abs{\ip f,\varphi_s.}^2\ge{}\frac{\sigma^2}4\abs{R_\T}, 
\end{equation}
$\abs{R_\T}$ is maximal, and ${\boldsymbol \omega}_\T$ is most clockwise. 
Then set $\tau(\T)\subset
\S^{\text{\rm stock}}$ to be the maximal (with respect to inclusion) tree in $\S^{\text{\rm stock}}$ with top 
${\boldsymbol \omega}_\T\times R_\T$.  Update 
\begin{equation*}
 \vZT_{+}=\vZT_{+}\cup\{\T\},\qquad \vZT_{\text{\rm big}}=\vZT_{\text{\rm big}}\cup\{\tau(\T)\},\qquad
 \S^{\text{\rm stock}}=\S^{\text{\rm stock}}-\T.
\end{equation*}
When  $\size {\S^{\text{\rm stock}}}<\sigma/2$, set $  \S^{\text{\rm stock}}=\vZT_{\text{\rm small}}$, and the process stops. 

To conclude the Lemma, we need to show that 
\begin{equation}\label{e.zs1}
\sum_{\T\in\vZT_{+}}\abs{ R_\T}\seq{}\sigma^{-2}.
\end{equation}
We have constructed these trees so that they satisfy 
a property very useful to the verification of this inequality.  Suppose $\T\not=\T'\in\vZT_{+}$ and $s\in\T$ 
and $s'\in\T'$ with $\bw
 {} \subset\bw {'1}$.  Then $R_\T\cap R_{s'}=\emptyset$.  We refer to this property as ``strong disjointness.''{} 
 To see that it holds, 
we have ${\boldsymbol \omega}_\T\subset \bw{'1}$, so that ${\boldsymbol \omega}_\T$ lies in the clockwise direction 
from  ${\boldsymbol \omega}_{\T'}$. Hence $\T$ was constructed before $\T'$. 
Thus, if $R_\T\cap R_{s'}\not=\emptyset$, we would conclude from \e.<property/ that $s'<\boldsymbol{\zw}_\T\times R_\T $. 
Hence,  $s'\in \tau(\T)$, and so $s'$ would have been removed from
$ \S^{\text{\rm stock}}$ and so could not be in $\T'$.

Adopt the notation $F(\S')=\sum_{s\in\S'}\ip f, \varphi_s .\varphi_s$.  
Now observe that for $\S_{+}=\bigcup_{\T\in \vZT_{+}}\T$, we have 
\begin{align*}
\sigma^2\sum_{\T\in\vZT_{+}} \abs{R_\T}\seq{}& 
	\sum_{s\in \S_{+}} \ip f,\varphi_s .\ip \varphi_s,f. 
	\\{}={}& \ip f, F(\S_{+}) .
	\\{}\le{}& \norm f.2.\norm F(\S_{+}).2.
\end{align*}
This follows by Cauchy--Schwarz.  Recall that $\norm f.2.=1$.   To conclude \e.zs1/ we need only show that 
\begin{equation}\label{e.zs2}
\norm F(\S_{+}) .2.\seq{}\sigma\Bigl[\sum_{\T\in\vZT_{+}}\abs{ R_\T}\Bigr]^{1/2}.
\end{equation}

\medskip

For $s\in\T$, set 
\md4
{\mathbf B}_=(s)=\{s'\in \S_{+}-\T\mid \bw{}=\bw '\},
\\
{\mathbf B}(s)=\{s'\in \S_{+}-\T\mid \bw{}\subset_{\not=}\bw '\}.
\emd
  Note that 
if $s'\in\S_+-\T$ is such that  $\ip \varphi_s,\varphi_{s'}.\not=0$, then $s'\in {\mathbf B}_=(s)\cup{\mathbf B}(s) $. 
And if $s'\not\in {\mathbf B}_=(s)\cup{\mathbf B}(s)$, then 
$\ip \varphi_s,\varphi_{s'}.=0$.  We estimate 
\md3
\lVert F(\S_{+})\rVert_2^2\seq{}&\sum_{\T\in\vZT_{+}} \lVert F(\T)\rVert_2^2 \label{e.top}
\\\label{e.middle}
&\quad{}+{} \sum_{s\in\S_{+}}\ip f,\varphi_s.\ip \varphi_s ,F({\mathbf B}_=(s)).
\\\label{e.bottom}
&\quad{}+{} \sum_{s\in\S_{+}}\ip f,\varphi_s.\ip \varphi_s ,F({\mathbf B}(s)).
\emd

It is a routine matter to verify that for each $1$--tree $\T\in\vZT_{+}$, 
\begin{equation*}
\norm F(\T) .2.\seq{} \sigma\abs{R_\T}^{1/2}.
\end{equation*}
Therefore, the right hand side of \e.top/ is no more than 
\begin{equation*}
\sum_{\T\in\vZT_{+}} \lVert F(\T)\rVert_2^2\seq{}\sum_{\T\in\vZT_{+}}\sigma^2\abs{R_\T}.
\end{equation*}

For the term \e.middle/, we use this  estimate for tiles $s,s'$ with $\bw{}=\bw'$,  
\md0
\abs{\ip \zvf_s,\zvf_{s'}.}\lesssim{}  \operatorname{Tran}_{c(R_s)}\operatorname{Dil}_{R_s}^0\zq (c(R_{s'})).
\emd
This with Cauchy--Schwarz gives 
\md2
\e.middle/\le{}&\Bigl[ \sum_{s\in\S_{+}}\abs{\ip f,\varphi_s.}^2 
	\sum_{s\in\S_{+}} \ABs{ \ip \varphi_s ,F({\mathbf B}_=(s)).}^2 \Bigr]^{1/2}
\\{}\lesssim{}& 
\sum_{s\in\S_{+}}\abs{\ip f,\varphi_s.}^2
\\{}\lesssim{}& \zs^2 \sum_{\T\in\vZT_{+}}\abs{ R_\T},
\emd
as required by \e.zs2/.

As for\e.bottom/, use Cauchy--Schwarz in $s$ to bound it by
 \begin{align*} 
 \sum_{s\in\S_{+}}\ip f,\varphi_s.\ip \varphi_s ,F({\mathbf B}(s)).
 {}\le{}& 
\Bigl[ \sum_{s\in\S_{+}}\abs{\ip f,\varphi_s.}^2 \sum_{s\in\S_{+}} 
	\abs{\ip \varphi_s ,F({\mathbf B}(s)).}^2 \Bigr]^{1/2}
\\{}\seq{}&\Bigl[ \sigma^2\sum_{\T\in\vZT_{+}}\abs{ R_\T} 
\sum_{s\in\S_{+}}\abs{\ip \varphi_s ,F({\mathbf B}(s)).}^2 \Bigr]^{1/2}
\end{align*}
We show that for each $\T\in\vZT_{+}$,
\begin{equation} \label{e.zs6}
\sum_{s\in\T}\abs{\ip \varphi_s ,F({\mathbf B}(s)).}^2\seq\sigma^2\abs{R_\T}
\end{equation}
This will complete the proof of \e.zs2/.

Observe that if we set ${\mathbf B}=\bigcup_{s\in\T}{\mathbf B}(s)$, we have 
\begin{equation*}
\ip \varphi_s,F({\mathbf B}(s)).=\ip \varphi_s,F({\mathbf B}). \,.
\end{equation*}
This is a consequence of the grid structure on the intervals $\boldsymbol \omega_s $.  
Moreover,all tiles $s\in\T\cup {\mathbf B}$ have $\bw{}\supset{\boldsymbol \omega}_\T$.  This 
has two implications.  The first is that   
 all rectangles $R_s$ can be regarded 
as rectangles with respect to a fixed set of coordinate axes, those for the top of the tree.
Let $\operatorname M$ denote the strong maximal function computed in 
these coordinate directions.

 The second, is that the strong disjointness 
condition applies to each pair of tiles in the collection $\mathbf B$. This yields the essential observation 
 that    the rectangles $\{R_s\mid s\in{\mathbf B}\}$ are pairwise disjoint, and do not intersect $R_\T$. 

Make a further diagonalization of the set ${\mathbf B}$.  Set ${\mathbf B}_1=\{s\in{\mathbf B}\mid R_s\subset 4R_\T\}$, and for 
$k>1$, set ${\mathbf B}_k=\{s\in{\mathbf B}\mid R_s\subset 4^k R_\T, \ R_s\not\subset 4^{k-1}R_\T\}$. Let us point out that 
\begin{equation*}
\norm F({\mathbf B}_k).2.\seq{}4^k \sigma \abs{R_\T}^{1/2}. 
\end{equation*}
Indeed,  recalling the notation \e.zqrs/ 
\begin{equation*}
\abs{F({\mathbf B}_k)}\seq{}\sigma\sum_{s\in{\mathbf B}_k} \zq_{R_s}*\ind R. 
\end{equation*}
Therefore, for any function $g$, 
\begin{align*}
\ip F({\mathbf B}_k),g.\seq{}&\sigma\int\sum_{s\in{\mathbf B}_k}\ind {R_s}\zq_{R_s}* g\; dx
\\{}\seq{}&\sigma\int \sum_{s\in{\mathbf B}_k}\ind {R_s} \operatorname Mg\; dx
\\{}\seq{}&\sigma4^k \abs{R_\T}^{1/2}\norm g.2.
\end{align*}
Clearly, this implies \e.zs6/ for ${\mathbf B}_1$.

 For $k>1$, we can strengthen our inequality to the following. 
\begin{equation*}
\norm F({\mathbf B}_k). L^2(4^{k-1}R_\T).\seq{} 4^{-10k}\abs{R_\T}^{1/2}.
\end{equation*}
We use this, together with the fact that the functions $\{ \zvf_s \mid s\in \T\} $ are equivalent to 
an orthonormal basic sequence.  Hence, 
\begin{equation*}
\sum_{s\in\T} \abs{\ip \varphi_s, F({\mathbf B}_k)\ind {4^{k-1}R_\T}.}^2\seq{}4^{-k}\abs{R_\T}. 
\end{equation*}
But the functions $\{\zvf_s\mid s\in\T\}$ are highly concentrated in a neighborhood of $R_\T$.  In particular, 
for any function $g$, 
\begin{equation*}
\sum_{s\in\T} \abs{\ip \varphi_s \ind {\{\ZR^2-4^{k-1}R_\T\}},g.}^2\seq{}4^{-10k}\lVert g\rVert_2^2. 
\end{equation*}
Clearly, this completes the proof of \e.zs6/.


\subsection*{Proof of \l.tree/}


\begin{figure}
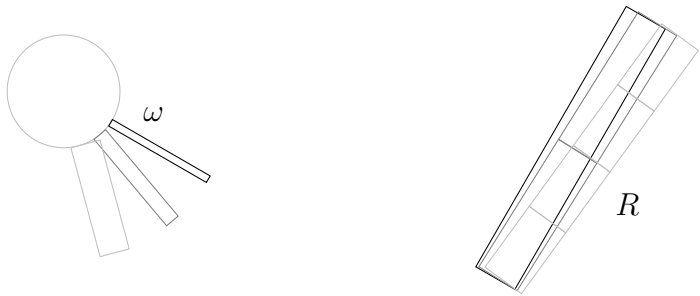
  \label{f.tree}
\begin{center} 
	\begin{pgfpicture}{0cm}{0cm}{5cm}{5cm}
	\begin{pgftranslate}{\pgfpoint{0cm}{2.5cm}}
	{\color{lightgray}
	\pgfcircle[stroke]{\pgfpoint{0cm}{0cm}}{0.75cm} 
	}
		\begin{pgfrotateby}{\pgfdegree{-30}}
		\pgfrect[stroke]{\pgfpoint{.75cm}{-.1cm}}{\pgfpoint{1.5cm}{0.1cm}}
		\end{pgfrotateby}
		\pgfputat{\pgfpoint{1.2cm}{-.3cm}}{\pgfbox[center,center]{$\omega$}}
		{\textcolor{gray}{ 
		\begin{pgfrotateby}{\pgfdegree{-50}}
		\pgfrect[stroke]{\pgfpoint{.75cm}{-.1cm}}{\pgfpoint{1.5cm}{0.2cm}}
		\end{pgfrotateby}
		}}
		{\textcolor{lightgray}{ 
		\begin{pgfrotateby}{\pgfdegree{-75}}
		\pgfrect[stroke]{\pgfpoint{.75cm}{-.1cm}}{\pgfpoint{1.5cm}{0.4cm}}
		\end{pgfrotateby}
		}}
	\end{pgftranslate}
	\begin{pgftranslate}{\pgfpoint{5.5cm}{-1cm}}
	\begin{pgfrotateby}{\pgfdegree{-30}}
	{
	\pgfrect[stroke]{\pgfpoint{-0.6cm}{1cm}}{\pgfpoint{0.6cm}{4cm}}
	}\end{pgfrotateby}
	{\textcolor{gray}{ 
	\begin{pgfrotateby}{\pgfdegree{-32}}
	\pgfrect[stroke]{\pgfpoint{-0.6cm}{1cm}}{\pgfpoint{0.6cm}{2cm}}
	\pgfrect[stroke]{\pgfpoint{-0.6cm}{3cm}}{\pgfpoint{0.6cm}{2cm}}
	\end{pgfrotateby}
	}}
	{\textcolor{lightgray}{ 
	\begin{pgfrotateby}{\pgfdegree{-36}}
	\pgfrect[stroke]{\pgfpoint{-0.6cm}{1cm}}{\pgfpoint{0.6cm}{1cm}}
	\pgfrect[stroke]{\pgfpoint{-0.6cm}{2cm}}{\pgfpoint{0.6cm}{1cm}}
	\pgfrect[stroke]{\pgfpoint{-0.6cm}{3cm}}{\pgfpoint{0.6cm}{1cm}}
	\pgfrect[stroke]{\pgfpoint{-0.6cm}{4cm}}{\pgfpoint{0.6cm}{1cm}}
	\end{pgfrotateby}
	}}
	\pgfputat{\pgfpoint{2cm}{2cm}}{\pgfbox[center,center]{$R$}}
	\end{pgftranslate}
	\end{pgfpicture}
	\end{center} 
	\caption{A few possible tiles in a $1 $--tree.  Rectangles $\omega_s $ are on the left in different shades of gray. 
	Possible locations of $R_s $ are in the same shade of gray. } 
	\end{figure}

We may fix the vector field, and assume that the standard basis $(e,e_\perp)$ are the basis for the rectangle $R_\T$.  As a consequence, we can without loss of generality assume that this is the basis for all the tiles $s\in\T$, and we write $R_s=R_{s,e}\times R_{s,e_\perp}$. 
Set $\zd=\dense \T $ and $\zs=\size \T $. 
For this proof, we set 
\md0
\zq^{(p)}=\operatorname{Tran}_{c(R_s)}\operatorname{Dil}^{p}_{R_s},\qquad 0<p<\infty. 
\emd

Let $\z0J_{e}$  be those maximal dyadic intervals $J$ in $\ZR$ for which for $3J$ does not contain an interval $R_{s,e}$
 for some $s\in\T$. This collection partitions \ZR. Let $\z0J_{e_\perp}$ be that partition of \ZR\ into maximal dyadic intervals $J$  such 
 that   $R_{\T,e_\perp}\not\subset3J$.  
Let $\z0K=\z0J_e\times \z0J_{e_\perp}$.  For each rectangle $K\in\z0K$, set $\T(K,\pm)$ to be those $s\in\T$ for which 
$\pm\log \abs{R_{s,e}}/\abs{K_e}>0$.

Choose signs $\zve_s\in\{\pm\}$ such that 
\md0
\sum_{s\in\T}\abs{\ipf \ip \zf_s,\ind F .}=\int_F\zve_s\sum_{s\in\T}\ipf \zf_s \; dx 
\emd
Let 
\md0
A^{\S'}=\sum_{s\in\S'}\zve_s\ipf \zf_s .
\emd

Estimate, without appeal to cancellation 
\md0
\int_K \abs{ A^{\T(K,-)}}\; dx
{}\seq{}
	\zs\int_K{\sum_{s\in\T(K,-)} \ind {\bw2}(v(x))   \zq^{(\zI)}_{R_s}(x)   }\; dx 
\emd
 The rectangles $R_s$ are smaller than those of $K$, and do not intersect it.   
For each $x \in K$, 
 the numbers   $ \zq^{(\zI)}_{R_s}(x)\seq{}1$, and they 
decrease as $\operatorname{dist}(R_s,K)\cdot\scl s$ increases. The integral above is at most 
\md2
\int_K \abs{ A^{\T(K,-)}}\; dx{}\seq{}\zd \zs\min(\abs{K}, \abs{R_\T} \sup_{s\in\T}\sup_{x\in K}\zq^{(\zI)}_{R_s}(x)).
\emd
This is summable over $K\in\z0K$ to at most $\seq{}\zd\zs\abs{R_\T}$.

\medskip 

If $\T(K,+)$ is non-empty, then $K\subset 4R_\T$.  Set 
\md0
G_K:=K\cap\bigcup_{s\in \T(K,+)}v^{-1}(\bw 2)
\emd
which contains the support of $\ind K A^{\T(K,+)}$.  Our assertion is that $\abs {G_K}\seq{}\zd\abs{K}$.  
To see this, let 
 $K=K_e\times K_{e\perp}$, and $K'_e$ be the dyadic interval that contains $K_e$ and $\abs{K'_e}=2\abs {K_e}$. 
This interval is small than that of $R_{\T,e}$.  
Then, by maximality, $3K'_{e}$ contains some $R_{s,e}$, for $s\in\T$.  Let $s'$ be the tile $s<s'<{\boldsymbol \zw}_\T
\times R_\T$ 
for which $R_{s'}=K'_e\times R_{\T,e_\perp}$.  We have $G_K\subset v^{-1}(\bw {'1})$.  And since $\dense \T\le\zd$, we 
deduce the important inequality that 
\md0
\abs{G_K}\le{} \abs{K\cap v^{-1}(\bw ')}\seq{}\zd\abs{K}.
\emd

If $T(K,+)$ is a $2$--tree, then the sets $\{ \bw 2\mid s\in\T(K,+)\}$ are pairwise disjoint, so that the set 
$v^{-1}(\bw 2)$ are either equal or disjoint.  Hence 
\md0
\ind K \abs{ A^{\T(K,+)}}\seq{}\zs\ind {G_K}.
\emd
Our desired bound $\int_K \abs{ A^{\T(K,+)}}\; dx\seq{}\zs\zd\abs K $ is immediate. 

If $\T(K,+)$ is a  $1$--tree, then all of the interval $\bw 2$ contain ${\z2 \zw}_\T$.  Thus, for each $x$, there are 
$\zve_{\pm}(x)$ so that $v(x)\in\bw 2$ iff $\zve_-(x)\le\scl s\le\zve_+(x)$.  In particular, if $v(x)\in\bw 2$, then we have 
$\abs{v(x)-e}\seq{}\zve_-(x)/{\mathsf {ann}} $.  This permits us to argue that the vector field $v(x)$ can be assumed 
to be constant.  Specifically one sees the inequality 
\md0
\ind {\bw 2}(v(x))\int_\ZR \abs{\varphi_{s}(x-yv(x))-\zvf_s(x-y e)}\abs{\zc_s(y)}\; dy\seq{} \frac{\zve_-(x)}{\scl s} \zq^{(2)}_{R_s}(x).
\emd
This inequality is in fact a  straight forward calculation, but one that depends very much on the uniform 
assumptions \e.zc-Space/.

We make more definitions, and in the first definition, replace the vector field by $e $.
\md2
\zf'_s(x):={}&\int \zvf_s(x-ye)\zc_s(y)\; dy,
\\
A_K:={}&\sum_{s\in\T(K,+)}\ip f,\zvf_s. \ind {\bw 2}(v(x)) \zf'_s, \qquad K\in\z0K , 
\\
 B:={}&\sum_{s\in\T}\ip f,\zvf_s.  \zf'_s.
\emd
In the last line, it is a matter of convenience to assume that all trees $\T(K,+) $ are \tree1s, and that 
$\T=\bigcup_{K\in\z0K} \T(K,+) $. 
We have, by a straightforward estimate, 
\md0
\Abs{A^{\T(K,+)}-A}\seq{}\zs\ind {G_K}
\emd
so that it suffices to estimate the $L^1$ norm of $A$ on  the interval $K$.  

There are two points: The first  is that for each $K $, $A_K$ is dominated by a maximal function applied to $B$. 
The second is that $B $ obeys the the inequality 
\md1  \label{e.B<}
\norm B.2. {}\lesssim{}\zs\abs{R_\T}^{1/2 }
\emd

\medskip 
Let us turn to the second claim.  The claim will follow immediately if we can see that the we have the universal 
estimate 
\md1 \label{e.BB} 
\Norm \sum_{s\in\T }a_s \zf_s' .2. {}\lesssim{} \Bigl[ \sum_{s\in\T} \abs{ a_s}^2\Bigr]^{1/2}. 
\emd
This holds for all numerical sequences $\{a_s\} $.

Indeed, the assumptions that we have placed on the functions $\varphi $, namely that 
it have Fourier support in a small ball around the origin, and the uniform assumption on the Fourier support of the 
functions $\psi_s $, as expressed in \e.zc-Fourier/,  imply that the functions $\zf_s' $ are supported in 
 $\zw_s $.\footnote{This calculation is available to us as we have replaced the vector field in the definition of 
$\phi_s' $ by a fixed vector $e $. } 
As these rectangles are either equal or disjoint, it suffices to prove the estimate \e.BB/ when the 
 tree $\T $ has $\zw_s=\zw_{s'} $ for all $s,\,s'\in\T $.    With this assumption, the assumption \e.zc-Space/ (with 
 the value of $\zk $ suitably small) 
 implies the estimate 
 \md0
 \abs{ \zf_s'} {}\lesssim{} \zq^{(2)}_{R_s} 
 \emd
And this is more than enough to prove \e.BB/ under this restrictive assumption on the tree. 

\medskip

Let us dominate $A_K $.  For any choice of $\zve_-<\zve_+$, we have 
\md0
\ind { \{\zve_-<\scl s<\zve_+\}} \zf'_s={}  (\zz_{\zve_-}-\zz_{\zve_+})*\zf'_s
\emd
in which we take $\zz$ to be  non--negative Schwartz function on the plane satisfying $\ind {[-1/2,1/2]^2}\le\widehat
\zz\le{} \ind{[-5/8,5/8]}$ and set 
$\zz_\zve(x_e,x_{e_\perp})=\zve \mathsf{ann} \zz(\zve x_e, \mathsf{ann} x_{e_\perp})$. The identity follows from 
the frequency properties of  $\zvf$ and of the class of functions $\zc_t $ as described in \e.zc-Fourier/.
From this, we conclude that 
\md0
\abs {A_K}\seq{} \sup_{J\supset R}\dashint_J \abs{B(z)}\; dz 
\emd
with the last quantity being a supremum over all rectangles $J $ that contain $K$. This supremum 
is constant on $K$.  Thus, 
\md0
\int_{G_K} \abs{A_K}\; dx\seq{} \abs{G_K} \inf_{x\in K}{\operatorname M B}(x)\seq{}\zd \abs{K}\inf_{x\in K}{\operatorname M B}(x)
\emd
where $M$ denotes the strong maximal function with respect to the $(e,e_\perp)$ coordinates.

This last estimate is to be summed  over $K\subset 4R_\T$. 
\md2
\zd\sum_{K\subset 4R_\T} \abs{K}\inf_{x\in K}{\operatorname M B}(x)
{}\lesssim{}& \zd \int_{4R_\T}\operatorname  MB \; dx 
\\{}\lesssim{}& \zd  \abs{R_\T}^{1/2} \norm \operatorname   MB .2. 
\\{}\lesssim{}& \zd \zs \abs{R_\T}
\emd
This completes the  the proof of the Tree Lemma.



\begin{bibsection}
 
 \begin{biblist}
 
 \bib{MR1955263}{article}{
    author={Barrionuevo, Jose},
    author={Lacey, Michael T.},
     title={A weak-type orthogonality principle},
   journal={Proc. Amer. Math. Soc.},
    volume={131},
      date={2003},
    number={6},
     pages={1763\ndash 1769 (electronic)},
      issn={0002-9939},
    review={1 955 263},
    eprint={math.CA/0201061 },
} 


 
\bib{bourgain}{article}{
    author={Bourgain, J.},
     title={A remark on the maximal function associated to an analytic
            vector field},
 booktitle={Analysis at Urbana, Vol.\ I (Urbana, IL, 1986--1987)},
    series={London Math. Soc. Lecture Note Ser.},
    volume={137},
     pages={111\ndash 132},
 publisher={Cambridge Univ. Press},
     place={Cambridge},
      date={1989},
    review={MR 90h:42028},
}

\bib{MR92g:42010}{article}{
    author={Bourgain, J.},
     title={Besicovitch type maximal operators and applications to Fourier
            analysis},
   journal={Geom. Funct. Anal.},
    volume={1},
      date={1991},
    number={2},
     pages={147\ndash 187},
      issn={1016-443X},
    review={MR 92g:42010},
}


 
  

\bib{carleson}{article}{
    author={Carleson, Lennart},
     title={On convergence and growth of partial sumas of Fourier series},
   journal={Acta Math.},
    volume={116},
      date={1966},
     pages={135\ndash 157},
    review={MR 33 \#7774},
}

\bib{christ}{article}{
    author={Christ, Michael},
     title={Personal Communication},
}  



 

\bib{fefferman}{article}{
    author={Fefferman, Charles},
     title={Pointwise convergence of Fourier series},
   journal={Ann. of Math. (2)},
    volume={98},
      date={1973},
     pages={551\ndash 571},
    review={MR 49 \#5676},
}

\bib{gl1}{article}{
    author={Grafakos, Loukas},
    author={Li, Xiaochun},
     title={Uniform bounds for the bilinear Hilbert transform, I},
   journal={Ann. of Math., to appear},
    }



\bib{gtt}{article}{
    author={Grafakos, Loukas},
    author={Tao, Terence},
    author={Terwilleger, Erin},
     title={$L\sp p$ bounds for a maximal dyadic sum operator},
   journal={Math. Z.},
    volume={246},
      date={2004},
    number={1-2},
     pages={321\ndash 337},
      issn={0025-5874},
    review={2 031 458},
} 

\bib{nets}{article}{
    author={Katz, Nets Hawk},
     title={Maximal operators over arbitrary sets of directions},
   journal={Duke Math. J.},
    volume={97},
      date={1999},
    number={1},
     pages={67\ndash 79},
      issn={0012-7094},
    review={MR 2000a:42036},
}


\bib{esi}{article}{
	author={Lacey, Michael},
	title={Carleson's Theorem: Proof, Complements, Variation},
	journal={Pub. Mat.},
	year={2004},
	}

\bib{laceyli}{article}{
	author={Lacey, Michael},
	author={Li, Xiaochun},
	title={On the Hilbert Transform on $C^{1+\ze}$ families of lines},
	journal={preprint},
	year={2003},
}

\bib{laceythiele}{article}{
    author={Lacey, Michael},
    author={Thiele, Christoph},
     title={A proof of boundedness of the Carleson operator},
   journal={Math. Res. Lett.},
    volume={7},
      date={2000},
    number={4},
     pages={361\ndash 370},
      issn={1073-2780},
    review={MR 2001m:42009},
}


\bib{MR81a:42027}{article}{
    author={Nagel, Alexander},
    author={Stein, Elias M.},
    author={Wainger, Stephen},
     title={Hilbert transforms and maximal functions related to variable
            curves},
 booktitle={Harmonic analysis in Euclidean spaces (Proc. Sympos. Pure Math.,
            Williams Coll., Williamstown, Mass., 1978), Part 1},
    series={Proc. Sympos. Pure Math., XXXV, Part},
     pages={95\ndash 98},
 publisher={Amer. Math. Soc.},
     place={Providence, R.I.},
      date={1979},
    review={MR 81a:42027},
}


 \bib{stein}{article}{
    author={Stein, E. M.},
     title={Problems in harmonic analysis related to curvature and
            oscillatory integrals},
 booktitle={Proceedings of the International Congress of Mathematicians,
            Vol. 1, 2 (Berkeley, Calif., 1986)},
     pages={196\ndash 221},
 publisher={Amer. Math. Soc.},
     place={Providence, RI},
      date={1987},
    review={MR 89d:42028},
}


\end{biblist}

\end{bibsection}
\end{document}